\input amstex
\documentstyle{amsppt}
\input bull-ppt
\keyedby{bull271e/jxs}

 \define\C{{\Bbb C}}
 \define\Z{{\Bbb Z}}
 
\define\Ger{\germ }
 \define\g{{\Ger g}}
 \define\h{{\Ger h}}
 
 \define\r{{\Ger r}}
 \define\s{{\Ger s}}
\define\germu{{\Ger u}}
\define\a{{\Ger a}}
 \define\so{{\Ger{so}}}
 \define\gsp{{\Ger{sp}}}
 \define\gsl{{\Ger{sl}}}
 \define\f{{\Ger f}}
 \define\rg{{R[\g]}}
 \define\ee{\varepsilon}
 \predefine\ge=\geq
 \redefine\ge{G^{\ee}}
 \define\omin {{O'_{\roman{min}}}}

 \define\too{\longrightarrow }
 \define\ov{\overline}
 \define\vs{\vskip .6pc}

\topmatter
\cvol{26}
\cvolyear{1992}
\cmonth{April}
\cyear{1992}
\cvolno{2}
\cpgs{269-275}

\title
 Nilpotent orbits, normality,\\
and Hamiltonian group actions\endtitle
\date March 29, 1991 and, in revised form, 
May 29, 1991\enddate
\author
Ranee Brylinski 
and
Bertram  Kostant 
\endauthor
\address Department of Mathematics, Pennsylvania State 
University, 
University Park, Pennsylvania 16802\endaddress
\address Department of Mathematics, Massachusetts 
Institute of Technology,
 Cambridge, Massachusetts 02139\endaddress
\thanks The first author is an Alfred P. Sloan 
Fellow\endthanks 
\thanks The work of the second author was supported in 
part by NSF Grant DMS-8703278\endthanks
\abstract Let  $M$ be a $G$-covering of a nilpotent orbit
in $\g$ where $G$ is a
complex semisimple Lie
group and $\g=\text{Lie}(G)$. We prove that
under Poisson bracket the space $R[2]$ of
homogeneous  functions on $M$
of degree 2  is the unique maximal
semisimple Lie subalgebra  of $R=R(M)$ containing $\g$.
The action of $\g'\simeq R[2]$ exponentiates
to an action of the
corresponding Lie group $G'$ on a
$G'$-cover
$M'$ of a nilpotent orbit
in $\g'$ such that  $M$ is open dense in $M'$.
We determine all such pairs $(\g\subset\g')$.\endabstract

\endtopmatter

\document

The theory of coadjoint orbits of Lie groups is
central to a number of areas in
mathematics.  A list of such areas would include (1) group
representation theory,  (2) symmetry-related Hamiltonian
mechanics and attendant physical theories, (3) symplectic 
geometry,
(4) moment maps, and (5) geometric quantization.
From many points of view
the most interesting cases arise when the group $G$ in 
question
is semisimple.
For semisimple $G$  the most familar
of the orbits are of semisimple elements.
In that case the associated representation theory is 
pretty much
understood (Borel-Weil-Bott and noncompact analogues, 
e.g., Zuckerman
functors).
The isotropy subgroups are
reductive  and the  orbits
are   in one form or another  related to flag
manifolds  and their natural generalizations.
 
A totally
different experience is encountered  with
nilpotent orbits
of semisimple groups. Here the associated representation
theory (unipotent representations) is poorly understood and
there is a loss of reductivity of isotropy subgroups. To
make matters worse (or really more interesting) orbits are 
no
longer closed and there can be a failure of normality for
orbit closures. In addition simple connectivity is
generally gone but more seriously there may exist no
invariant polarizations.
 
The interest in   nilpotent orbits of semisimple Lie groups
has increased sharply over the last two decades.
This perhaps  may be attributed to  the reoccuring 
experience
that sophisticated aspects of semisimple group theory
often leads one to these orbits  (e.g., the Springer
correspondence with representations of the Weyl group).
 
In this note we announce new results concerning the
symplectic and algebraic geometry of the nilpotent orbits 
$O$
and the symmetry groups of that geometry.
The starting point is the
recognition (made also by others)  that the ring $R$ of
regular functions on  any $G$-cover $M$  of $O$ is not 
only a
Poisson algebra (the case for any coadjoint orbit) but 
that $R$
is also naturally graded.  The key theme  is that the same
nilpotent orbit may be \lq\lq shared" by more than one 
simple
group, and the key result is the determination of all pairs
of simple Lie groups having a shared
nilpotent orbit.  Furthermore  there is  then a
unique maximal such group and this group is encoded in the
symplectic and algebraic geometry of the orbit.  Remarkably
a covering of nilpotent orbit of a classical
group may \lq \lq see" an exceptional Lie group as the
maximal symmetry group of this sympectic manifold.  A
beautiful instance of this is  that $G_2$ is the
symmetry group of the simply connected covering  of the
maximal nilpotent orbit of $\operatorname{SL}(3,\Bbb C)$ 
and that this
six-dimensional space \lq \lq becomes" the minimal nilpotent
orbit of $G_2$.
 
Our work  began  with a desire to thoroughly investigate
a striking discovery of Levasseur,
Smith, and Vogan. They found that the failure of the closure
of the eight-dimensional nilpotent orbit of $G_2$ to be a 
normal
variety may be \lq \lq remedied" by refinding this orbit as
the minimal nilpotent orbit of $\operatorname{SO}(7,\C)$. 
The failure has a lot
to do with the seven-dimensional representation of $G_2$. In
general  given $M$
we have found that there exists a  unique  minimal
representation $\pi$  (containing the adjoint) wherein $M$
may be embedded with normal closure.  It was the study of
$\pi$ that led to the discovery of the maximal symmetry
group $G'$. Using a new general transitivity result
for coadjoint orbits  we prove
that, modulo a possible normal Heisenberg subgroup (and that
occurs in only one case), $G'$ is semisimple.
 
Past experience has shown that the action of a subgroup
$H$  on a coadjoint orbit of $G$ is a strong
prognostigator as to how the corresponding representation
$L$
of $G$ decomposes under $H$. If this continues to hold for
unipotent representations our classification result should
yield all cases where $L$  remains
irreducible (or decomposes finitely) under a semisimple
subgroup.

\heading
1.  The maximal symmetry group and
``shared" orbits \endheading
Let $G$ be a simply connected complex semisimple Lie group
and $\g$ the Lie algebra of $G$.
  Let  $e\in\g$ be nilpotent and assume
(for simplicity of exposition but with no real loss) that
$e$ has nonzero projection in every simple component of 
$\g$.
Let $O$ be the adjoint orbit of $e$ and let
$\nu\:M\too O$
be a $G$-covering.
Let  $R=R(M)$  be the ring of regular
functions on $M$  and let
$\rg$ be the copy of $\g$ in $R$  defining the
moment map  $M\too \g^*$.
Identify  $\g^*\simeq\g$
in a $G$-equivariant way.
 
 Now $R$ carries a   $G$-invariant ring grading
 $R=\bigoplus_{k\geq 0}R[k]$
($k\in\Z$) such that $\rg\subset R[2]$.
 Then $R[0]=\C$.   The Poisson bracket satisfies
$[R[k],R[l]]\subset R[k+l-2]$,  for all $k,l$ (see also 
[7]).
 Hence  $R[2]$ is a
 finite-dimensional
 Lie subalgebra of $R$ under Poisson bracket.
 Our first main result is

\proclaim{Theorem 1} 
\roster\item"(i)"
\<$R[2]$ is  a semisimple Lie algebra, call it $\g'$.
If $\g$  is simple then  $\g'$ is  simple .
\item"(ii)"  The condition
$\rg\subset R[2]$ determines the
$G$-invariant ring grading on $R$ uniquely.
\item"(iii)" \<$R[2]+R[1]+R[0]$ is the unique maximal
finite-dimensional Lie subalgebra of $R$
containing $\rg$.
\endroster

\endproclaim


\midinsert
\toptablecaption{\smc Table 1}
$$\vcenter{ \settabs\+
 aaaaaaaaaa& aaaaaaaaaaaaaa& aaaaaaaaaaaaaaa&
aaaaaaaaaaa\cr
 \+  & $\g$ & $\g'$ & V\cr
 \vskip .5pc \hrule \vskip .5pc
 \+ (1)& $G_2$   & $\so(7)$ & $\C^7$  \cr
 \+ (2)& $\so(2n+1)$& $\so(2n+2)$&  $\C^{2n+1}$    \cr
 \+ (3)& $\gsp(2n)$& $\gsl(2n)$& $(\wedge^2\C^{2n})/\C$ \cr
 \+ (4)& $F_4$       & $E_6$     &  $\C^{26}$   \cr
 \+ (5)& $\gsl(3)$ & $G_2$  &  $\C^3\oplus\wedge^2\C^3$ \cr
 \+ (6)& $\so(2n)$& $\so(2n+1)$&  $\C^{2n}$  \cr
 \+ (7)& $\so(9)$ & $F_4$    &  $\C^{16}$  \cr
 \+ (8)& $\so(8)$ & $F_4$ & $\C^8\oplus\C^8\oplus\C^8$  \cr
 \+ (9)& $G_2$  & $\so(8)$ & $2\C^7$   \cr}
$$ \vs
\endinsert
In Table 1 all the Lie algebras are complex simple; 
particularly, in 
(2) and (3) $n\geq 2$ and in (6) $n\geq 3$.
 The last column $V$ is a representation of $\g$, written 
as a
sum of its irreducible components (only
 fundamental dominant representations occur here).
 In (3) $V\oplus L=\bigwedge^2\C^{2n}$
 where  $L\simeq\C$ is  $\gsp(2n)$-invariant.
 In (8) $V$ is the sum of the
 standard and the two half-spin representations.
 
\proclaim{Theorem 2} 
Table  {\rm 1} gives a complete
 list of the simple Lie algebra
 pairs $(\g\subset\g')$
 that arise in Theorem \RM1 with $\g\neq\g'.$
Furthermore
 $\g'=\g\oplus V$
 as $\g$-modules.
\endproclaim

Now the Poisson bracket on $R$ defines an alternating
bilinear form  $\beta$ on
$R[1]$  and  a Lie algebra homomorphism
$R[2]\to\gsp\,\beta$.

\proclaim{Theorem 3} 
$\beta$ is a symplectic form on
$R[1]$  so that
$R[1]\oplus R[0]$ is a Heisenberg Lie algebra.
If ${\g}$ is simple
then $R[1] \ne 0$ in one and only one case, namely, when
${\g}$ is of type $C_n$
for some $n$ and
$M$ is the simply connected \RM(double covering\/\RM) of
 the minimal nontrivial nilpotent
orbit of ${\g}$. In
that case
$R[2]\simeq\gsp(2n)$ and
$R[1]\simeq\C^{2n}$
generates $R$ freely.
\endproclaim

Now the functions in $\g'\simeq R[2]$
 define a map
 $\phi\:M\too\g'$
 (again identify $(\g')^*\simeq\g'$ ).
 Let $G'$ be the simply connected
 Lie group with Lie algebra $\g'$.
The next result says that up to
birationality
{\it $\phi$  is  a moment map for $G'$}.
 
\proclaim{Theorem 4}  The image  $\phi(M)$
 lies in a nilpotent orbit
 $O'$ of $G'$
 and
 $\phi(M)$
 is  Zariski open dense in $O'$.
 There
 exists a unique
 $G'$-covering
 $\nu'\:M'\too O'$
 such that $M'$ contains $M$
 and $\nu'$ extends $\phi$.
 Moreover   $M$ and $M'$
 have the same
 regular functions, that is,
 $R(M)=R(M')$, and also the
 same fundamental groups,
 that is,
 $\pi_1(M)\simeq\pi_1(M')$.\endproclaim

We  can construct $M'$ in the following way.
Given $M$,  let $X=\operatorname{Spec} R$
 be the maximal ideal spectrum of the finitely generated
$\C$-algebra $R$. Then $X$ is a normal affine variety.
Furthermore $X$ contains  $M$
as an open dense subset. We call
$X$ the {\it normal closure\/} of $M$.
Indeed if $M=O$ and
$\ov O$ is normal then $X=\ov O$.

 Our construction is:
 {\it $G'$ acts on $X$
 and  $M'$ is the unique
 Zariski
 open orbit of $G'$ on $X$.}
Note that
$M$ is the unique Zariski
open orbit of $G$ on $X$ so that clearly
$M\subset M'\subset X$.

Thus the pair
$(M\subset M')$  constitutes an orbit (cover)
\lq\lq shared"  between $G$ and  $G'$ .
Moreover $M$ and $M'$ have the same normal closure
so that $X$ is  exactly shared between $G$ and  $G'$.
 
Now even though $X$ may be singular we will say
an isomorphism of $X$ is symplectic if the
corresponding automorphism of $R$ preserves the Poisson
bracket structure.

\proclaim{Theorem 5}  Assume  $R[1]=0$.
Then $X$ is a singular variety.
 The action of any connected
 Lie group of holomorphic
 symplectic isomorphisms of $X$
 that  extends the action of
 $G$ is given by a subgroup
 of $G'$.\endproclaim

Hence, assuming $\g$ is simple,
{\it the example of Theorem \RM3 is the one and
only one choice of $M$ such that $X$ is smooth,
and in that case  one has  $X\simeq\C^{2n}$}.

\heading 2.   Explanation of the table\endheading
 We now
 describe for each of the nine cases in Table 1 
 a choice of $M$ such that
 $(M,\g)$ gives rise to $(M',\g')$.
In each case $M'=O'$ is the
orbit of the highest
root vector in the simple Lie algebra
$\g'$, this is the minimal nontrivial
nilpotent orbit (cf. the next section).
 
In cases (1)--(4) $\g$ is any one of the four
simple  Lie algebras that are
doubly laced (i.e., having two root lengths).
Choose $O$ to be the orbit of a short root vector and let
$M$ be its simply connected cover. Then 
$M=O$ in (1) while $M$ is a  two-fold cover of $O$ in
(2)--(4).  Then
$\g'$ is simply laced .
 Furthermore
 $V=V_{\alpha}$ is the irreducible $\g$-representation
with highest weight $\alpha$ equal to the highest
short root of  $\g$.
Case  (1) is  a
restatement in our language of a  result
 proved by
 Levasseur and Smith [4] in answer to a conjecture of
 Vogan [6] (see introduction).
 
In  (5) choose  $O$  to be the six-dimensional
 (maximal) nilpotent
 orbit of all principal nilpotent
 elements and  let
 $M$ be the three-fold
 simply connected covering
 space of $O$.
This  case was
 discovered by us in
 collaboration with Vogan.
A noncommutative analog of this example
is given in a result of Zahid [8] .
 
 In (6)  choose $O$ to be the nilpotent orbit of Jordan type
(see e.g., [3])
$(3,1^{2n-3})$  and let $M$ be the
simply connected double cover of $O$.
In (7)  choose $O$ to be of Jordan type
$(2^4,1)$  and let $M$ be the
simply connected (double cover) of $O$.
In (8)  choose $O$ to be of Jordan type
$(3,2^2,1)$  and let $M$ be the
simply connected  (four-fold) cover of $O$.
 
In (9)
 choose  $O$ to be the unique ten-dimensional
 nilpotent orbit  and let
$M$ be the
simply connected  six-fold cover of $O$
$(\pi_1(O) \simeq S_3)$.
In this example
 Levasseur and Smith  already showed in [4] that
 $G$  has an open dense orbit on $O'$,
 again in response to a question of Vogan.
 Moreover Vogan  has constructed a
 unipotent representation theoretic analogue of the 
example. If
 $\pi$ denotes the minimal unitary representation of
 $\operatorname{SO}(4,4)$ (see e.g., [2]) then $\pi$ 
extends to the
 outer automorphism group $A$ of $\g$ and in particular to a
 group $S\simeq S_3$ that induces $A$.
 Vogan shows that a split form $G_o$ of $G_2$
 and $S$ behave like a Howe pair with respect to $\pi$ and 
that
 $\pi|G_o$ decomposes into
 six irreducible components.
 Furthermore McGovern in [5, Theorem 4.1] has constructed a
 Dixmier algebra analogue of this example

 Remarkably three of the four nilpotent orbits of $G_2$ have
now appeared as \lq\lq shared" orbits
(the principal orbit does not appear).
 
Cases (5)--(8)
are precisely the  the pairs
 $(\g\subset\g')$  in the table with
 $\hbox{rank}\,\,\g$
 =$\hbox{rank}\,\,\g'$.
Each pair  $(\g\subset\g')$ is of the form
$(\s_0\subset\s)$
where $\s$\ is a doubly laced simple Lie algebra
and $\s_0$ is a subalgebra of $\s$
containing a Cartan subalgebra of
$\s$ and all associated long root vectors.
Moreover every such pair
$(\s_0\subset\s)$
where  $\s_0$ is simple arises in (5)--(8).
On the other hand the
pairs $(\s_0\subset\s)$
where  $\s_0$ is nonsimple occur
precisely when
$\s\simeq\gsp(2n)$
and
$\s_0\simeq\gsp(2n_1)\oplus
 \cdots\oplus\gsp(2n_k)$
 where
 $n_1+\cdots+n_k = n$  and $k\geq 2$.
These pairs do arise in Theorem 1
(when $M$ is chosen so that in each simple
component one has the example of
Theorem 3) and  these together with (5)--(8)
exhaust all equal rank pairs
$(\g\subset\g')$  arising in Theorem 1 such that
$\g'$ is simple.
 
A general result regarding the ranks of
$\g$ and $\g'$ is that
$\operatorname{rank}\g'>\operatorname{rank}\g$ whenever 
$M=O$
and  $\g\neq\g'$.

Two instances of a triple of Lie algebras
having a
\lq\lq shared" orbit  can be found among these
 examples, namely,
(a) $\so(8)\subset\so(9)\subset\f_4$
(where $\f_4$ is of type $F_4$)
and (b)
$\g_2\subset\so(7)\subset\so(8)$.
This is  not unexpected by
the theory since we in fact prove that
 if $\h$ is $any$ Lie subalgebra between
$\g$ and $\g'$ then $\h$ is semisimple
and also
$\h$ is simple if $\g$ is simple.
Moreover if $H$ is the simply connected
group corresponding to $\h$
then the statements made in Theorem 4
and immediately afterward
for $G'$  and $X$ apply equally well to $H$ and $X$.
In particular one has a unique
open $H$-orbit
$M^{\h}$ in $X$,
$M\subset M^{\h}\subset M'$,
and
$M^{\h}$ covers $H$-equivariantly
a nilpotent $H$-orbit
$O^{\h}\subset\h$. Furthermore
the whole graded Poisson ring structure on $R$
arising from
$M$ and $O$ is the same as would arise from
$M^{\h}$ and $O^{\h}$.
In particular the maximal semisimple
Lie algebra $R[2]$ remains the same.

\heading 3. Methods of proof\endheading
Two key ideas are used in  proving the classification of 
pairs.
The first is  a representation theoretic.
We are able to
compute $\g'\simeq R[2]$ as a $\g$-module (for arbitrary 
$\g$ and $M$).
Let $\ee$ be a point of $M$ lying over $e$
and let $(h,e,f\,)$ be a standard basis of an  $\gsl(2)$ 
subalgebra
$\a$ of
$\g$.  For any $\g$-module  $V$  let
$V[2]$ be the  2-eigenspace of $h$ in the
fixed space $(V^*)^{G^\ee}$.
Then
$\g'\simeq\g\oplus n_iV_1\oplus\cdots\oplus n_sV_s$
where $V_1,\dots,V_s$  is a complete list of inequivalent 
simple
$\g$-modules, excluding  components of the adjoint
representation, such that
$n_i=\dim V_i[2]$ is nonzero.
This follows by recognizing that the grading on $R$
comes from exponentiating a natural action of a Cartan
subalgebra in  $\a$ .
 
The second
idea is due to David Vogan who observed that for a given
pair $(\g,\g')$  arising in Theorem 1
we may change (if necessary) the choice of
$O$ and $M$ so that $O'$  is minimal .  Vogan himself has
determined the pairs $(\g,\g')$ in many of the cases
listed  above.
 
A principle used in setting up the theory
is that one should study minimal embeddings of
$X$ in order to study $M$ (again $M$ is arbitary).
We prove the  covering map $\nu$ extends
to a finite $G$-morphism
$\overline {\nu}\:X \to \overline O$
and then the fiber
$\nu^{-1}(0)$ over zero is a single point $o$.
Then $o$ is the unique $G$-fixed point in $X$  and
the  maximal ideal of $R$
corresponding to $o$ is
$m = \bigoplus_{k=1}^{\infty}R[k]$.
Furthermore,  $X$ is singular if and only if
$X$ is singular at $o$ and then
$o$ is  the ``most" singular point of $X$
(cf. Theorems 3 and 5).
 
We say that
{\it the pair $(v,V)$ defines an embedding of $X$}
in case $V$ is a $G$-module and $v \in V^{G^{\ee}}$
(where $\ee\in M$)
 is such
that the natural map $M \to G\cdot v$ extends to a
$G$-isomorphism $X \to \overline {G \cdot v} $.
We  find that the Zariski
tangent space $T_o(X) = (m/m^{2})^*$
at  $o$
 provides a minimal embedding for $X$.

\proclaim{ Theorem  6}  There exists a vector $u \in
T_o(X)^{G^{\ee}}$ such that {\rm(i)} $(u,T_o(X))$ defines an
embedding of $X$ and {\rm (ii)} if $(v,V)$ is any pair
that defines an embedding of $X$ then there exists a
surjective $G$-map $\tau\:V \to T_o(X)$ such that
$\tau(v) = u$.\endproclaim

Clearly one has an injection
$R[2]\to m/m^2 $ when  $R[1]=0$.
Regarding  the normality of  $\ov O$ we find
that
{\it
$\overline O$ is a normal variety if and only if  as
$G$-modules $\rg \simeq  R[2]\simeq m/m^2 $.}

The proofs of Theorems 1, 3, and 5 are all applications
of the  following general
transitivity theorem  for coadjoint orbits
of a special kind of Lie algebra.
This  is one of our main results.

\proclaim{ Theorem 7} Assume that $\s$ is a 
finite-dimensional Lie algebra
over  ${\Bbb R}$ or ${\Bbb C}$
and that ${\s}$ is a semidirect sum
${\s = \r + \germu}$ where ${\germu}$ is an abelian
ideal in ${\s}$,  ${\r}$ is a semisimple
Lie subalgebra of ${\s}$, and   ${\germu}^{\r}=0$.
One may regard $\s^*=\r^*+\germu^*$ in an obvious way.
Let $\gamma\in\s^*$ and write
$\gamma=\mu +\lambda$ with $\mu \in\r^*$ and 
$\lambda\in\germu^*$.
Then one has
${\s}\cdot \gamma = {\r} \cdot\gamma$
if and only if $\lambda = 0 $.\endproclaim

Theorem 7 says that if
$\lambda \ne 0$ then the subgroup of
$\operatorname{Ad}{\s}$ corresponding to ${\r}$ cannot
operate transitively (even infinitesmally) at
$\gamma$ on
the coadjoint orbit.

\heading 4.  An application to symmetry of flag 
varieties\endheading
Finally we give an application of our results to
a well-known problem in geometry.
If $P$ is a
parabolic Lie subgroup of a simple Lie group $G$ it is a
solved problem (see [1]) to determine the connected
component of the full group $F$ of holomorphic
automorphisms of the projective variety $G/P$. It is
precisely for the $G$ given in
cases (1), (2), and (3) in the table that $P$ exists
so that $F$ is larger than that given by the action of
$G$.  Furthermore in those cases $F$ is in fact given by
the action of $G'$.
We   obtain  a
stronger statement  (and recover the known result)
in

\proclaim{Theorem 8}  Let $P$ be any parabolic subgroup of 
$G$ and choose $M$ to
be the unique open orbit in $T^*(G/P)$  so that we can
take $O$ to be the $G$-orbit in $\g$ of a Richardson
element in the nilradical of $\operatorname{Lie} P$. Then 
one has a
desingularization map $T^*(G/P) \to X$ and the pullback of
$R$ is the full  ring of regular functions on $T^*(G/P)$.
Furthermore the action of $G'$ on $X$ lifts uniquely to an
action as a group of symplectic holomorphic automorphisms
of $T^*(G/P)$ and as such  $G'$ is maximal.   $G'$
preserves the cotangent polarization of $T^*(G/P)$ so that
$G'$ acts on $G/P$ and hence there exists a parabolic
subgroup $P' \subset G'$ such that $G/P = G'/P'$.  The 
action of $G'$ on $G/P$
is the connected component of the group of all holomorphic
automorphisms of $G/P$.
Consequently
any connected Lie group of
symplectic holomorphic automorphisms of $T^*(G/P)$
containing the action of $G$ automatically preserves the
cotangent space polarization of $T^*(G/P)$ and
consequently will act as a group of holomorphic
automorphisms of $G/P$.\endproclaim

\heading Acknowledgments\endheading  

The authors thank Madhav Nori and
David Vogan for helpful discussions. The first author 
thanks the
Laboratoire Math\'ematiques Fondamentales of Universit\'e
Paris VI
for  its hospitality.

\enddocument